\documentclass[12pt,letterpaper]{article}
\usepackage[margin=1in]{geometry}
\usepackage[T1]{fontenc}
\usepackage{mathpazo}
\usepackage{setspace}
\usepackage{microtype}
\usepackage{needspace}
\usepackage{xcolor}
\definecolor{revisionblue}{RGB}{0,0,139}

\usepackage{amsmath,amsthm, amssymb,tikz,enumitem}
\usetikzlibrary{arrows.meta,decorations.pathreplacing,calc}
\usepackage[authoryear,round,longnamesfirst]{natbib}
\usepackage[hidelinks]{hyperref}
\expandafter\def\csname biburl@B13\endcsname{https://arxiv.org/abs/1109.5668}
\expandafter\def\csname biburl@BCK14\endcsname{https://arxiv.org/abs/1205.6721}
\expandafter\def\csname biburl@B16\endcsname{https://arxiv.org/abs/1406.5660}
\expandafter\def\csname biburl@BK18\endcsname{https://arxiv.org/abs/1708.02134}
\expandafter\def\csname biburl@BT19\endcsname{https://arxiv.org/abs/1805.11989}
\expandafter\def\csname biburl@EKMS00\endcsname{https://annals.math.princeton.edu/articles/11101}
\expandafter\def\csname biburl@JM02\endcsname{https://doi.org/10.1016/S0304-4149(01)00142-9}
\expandafter\def\csname biburl@JM04\endcsname{https://arxiv.org/abs/math/0301055}
\expandafter\def\csname biburl@hambly-martin\endcsname{https://arxiv.org/abs/math/0604189}
\expandafter\def\csname biburl@CGK93\endcsname{https://doi.org/10.1214/aoap/1177005277}
\expandafter\def\csname biburl@GK94\endcsname{https://doi.org/10.1214/aoap/1177005201}
\expandafter\def\csname biburl@ganguly2024last\endcsname{https://arxiv.org/abs/2411.08018}
\expandafter\def\csname biburl@LP18\endcsname{https://doi.org/10.1017/9781316104477}

\usepackage{mathtools}

\hypersetup{pdftitle={The shape function of last passage percolation up to constants}}
\hypersetup{pdfauthor={Jnaneshwar Baslingker, Daniel Spivak, B\'alint Vir\'ag}}

\newcommand{\E}{\mathbb{E}}
\newcommand{\R}{\mathbb{R}}
\newcommand{\Z}{\mathbb{Z}}
\renewcommand{\P}{\mathbb{P}}
\newcommand{\dd}{\,\mathrm d}
\def\l{\left}
\def\r{\right}

\newtheorem{theorem}{Theorem}
\newtheorem{corollary}[theorem]{Corollary}
\newtheorem{lemma}[theorem]{Lemma}
\newtheorem{proposition}[theorem]{Proposition}

\theoremstyle{remark}
\newtheorem{remark}[theorem]{\bf Remark}

\begin{document}
\title{{The shape function of last passage percolation\\up to constants}}
\author{Jnaneshwar Baslingker \and Daniel Spivak \and B\'alint Vir\'ag}
\date{}

\maketitle
\begin{abstract}
% We determine the diagonal shape function of lattice last passage
% percolation with nonnegative i.i.d.\ weights up to dimension-dependent
% constants, giving an exact tail-integral criterion for finiteness.
% For critical tails of exponent $q$ in dimension $q\geq2$, the mean
% passage time has order $n(\log n)^{(q+1)/(2q)}$. This proves sharpness
% of the logarithmic exponent obtained in the lower bound of Ganguly,
% Ginsburg and Nam, removes their iterated-logarithmic loss, and provides
% a matching upper bound. The proof uses entropy-controlled Poisson last
% passage percolation, whose quadratic penalty makes the curvature
% dependence explicit, and transfers the resulting estimates to the lattice.
We determine the diagonal shape function of lattice last passage percolation up to dimension-dependent constants and characterize its finiteness by a tail integral. The proof transfers estimates from an entropy-controlled Poisson model whose quadratic penalty makes the curvature dependence explicit. For critical tails of exponent \(q\) in dimension \(q\ge2\), the mean passage time has order \(n(\log n)^{(q+1)/(2q)}\), proving sharpness of the logarithmic exponent in the lower bound of Ganguly, Ginsburg and Nam.\end{abstract}

\section{Introduction}

For a positive Borel measure $\nu$ on $[0,\infty)$ and $q\geq2$, define the rate
\begin{equation}\label{eq: gen dim rate}
    J_q(\nu)
=
\int_0^\infty
\nu(x,\infty)^{2/(q+1)}x^{(q-1)/(q+1)}\,dx.
\end{equation}

Consider up-right paths in $\Z^2$ from $(1,1)$ to $(n,n)$, with i.i.d.\
nonnegative vertex weights of law $\nu$. Let $L_{\Z^2}(\nu,n)$ be the
maximum collected weight, and write
\[
\mu_{\Z^2}(\nu)
=
\lim_{n\to\infty}\frac{\E L_{\Z^2}(\nu,n)}{n}
\in[0,\infty].
\]

\begin{theorem}\label{thm:bounds for lattice model}
There exist universal constants $c,\bar c>0$ such that, for every
probability measure $\nu$ on $[0,\infty)$,
\[
c\,J_2(\nu)^{3/4}
\leq\mu_{\Z^2}(\nu)
\leq\bar c\,J_2(\nu)^{3/4}.
\]
\end{theorem}

Thus $J_2(\nu)<\infty$ is the exact condition for finiteness of the
diagonal shape function. The same quantity determines its size up to
constants uniform over the weight distribution. {The shape limits exist in $[0,\infty]$ by superadditivity of the expected passage values.}

\subsection{Curvature and the proof}

The main idea is to work first in an entropy-controlled model, where
the dependence on curvature is explicit. A path collects weighted
Poisson points and pays the quadratic penalty
\[
d\int_0^t |f'(u)|^2\,du.
\]
A constant-slope shear preserves the law of the environment. Changing
the terminal slope from zero to $v$ therefore changes the passage
value in distribution by exactly $-dtv^2$. The curvature is directly
controlled by $d$, rather than hidden in an unknown lattice shape
function. This explicit dependence is what allows us to combine
contributions from different weight scales without losing the sharp order.

For the upper bound, we expose the weight levels from largest to
smallest. An exponential insertion estimate bounds the contribution
of each new level, uniformly over the path through the previously
exposed points. The resulting entropy-controlled upper bound is then
transferred to the lattice by a coupling. For the lower bound, we redo
the recursive construction in the ordinary compound Poisson setting:
smaller-weight points are inserted between previously selected points,
with constraints on the changes in slope. These constraints preserve
enough space for subsequent insertions (see Figure~\ref{fig:multiscale}). We then transfer the Poisson
lower bound to the lattice. 
\subsection{Related work}

Early general bounds came from the theory of greedy lattice animals.
\citet{CGK93} proved linear upper bounds under a moment condition which,
in dimension two, is
\[
\int_1^\infty x(\log x)^{2+\varepsilon}\nu(x,\infty)\,dx<\infty
\qquad\text{for some }\varepsilon>0.
\]
\citet{GK94} established laws of large numbers for greedy lattice animals
and self-avoiding paths under the same condition. Since an up-right
path is a lattice animal, these estimates also give sufficient
conditions for finite last-passage growth.

\citet{JM02,JM04} weakened the sufficient condition to
\begin{equation}\label{eq: prev suff condition}
    \int_0^\infty\nu(x,\infty)^{1/2}\,dx<\infty
\end{equation}
and studied the resulting shape function, including its continuity
at the boundary. A finite second moment is necessary, leaving a gap
between the necessary and sufficient conditions. Theorem~\ref{thm:bounds for lattice model} identifies the exact
condition between these two assumptions.

\citet[Theorem~5]{ganguly2024last} showed that a finite second moment is
not sufficient, answering a question of \citet{JM02}. Their counterexamples
have tails
\[
\nu(x,\infty)\sim\frac{c_0}{x^2(\log x)^\beta},
\qquad c_0>0,\quad 1<\beta<\frac32.
\]
These laws have finite second moment, but
$L_{\Z^2}(\nu,n)/n\to\infty$ almost surely. Their quantitative lower bound
\citep[Proposition~4.7]{ganguly2024last}, valid with probability tending to
one, already has the correct logarithmic exponent: it is within a factor
$\log\log n$ of our lower bound for the mean, of order
$n(\log n)^{3/4-\beta/2}$. Their counterexample range is also sharp up to
the endpoint: our criterion gives finiteness exactly when $\beta>3/2$,
and an infinite shape function at $\beta=3/2$.

The continuum model used in our proof has its origins in quadratic
random-action problems associated with randomly forced Burgers
equations. \citet{EKMS00} developed the variational approach to
stationary solutions in the spatially periodic setting, using
one-sided action minimizers. \citet{B13} studied Poisson forcing on
the real line under a spatial concentration condition.
\citet{BCK14} then treated a space-time homogeneous Poisson environment,
where paths minimize a quadratic kinetic cost minus the number of
visited Poisson points. Up to normalization, this is the unit-weight
counterpart of the penalized model used here. They also exploited
shear invariance to obtain an exactly quadratic dependence of the
shape function on slope. \citet{CP11} proved shape theorem under the assumption \eqref{eq: prev suff condition}, for rate-one planar Poisson process with i.i.d. weights. \citet{B16} developed the corresponding
theory for random kick forcing, and \citet{BK18} discussed its
connections with random Hamilton--Jacobi equations and KPZ universality.

\citet{BT19} introduced the terminology \emph{entropy-controlled
last-passage percolation}, with applications to directed polymers with
heavy-tailed disorder. We adopt their terminology for the penalized model,
while calling the quadratic path functional its Dirichlet energy.

\subsection{Entropy-controlled last-passage percolation}

For the following two continuum models, let $\nu$ be a locally finite
positive Borel measure on $(0,\infty)$.
Let $\mathcal P$ be a Poisson process of marked points $(z,u,w)$ with
intensity
\[
dz\,du\,\nu(dw).
\]
For finite $\nu$, this is a compound Poisson environment with point
intensity $|\nu|$ and mark law $\nu/|\nu|$. For $d>0$, define
\[
L_D(\nu,d,(y,t))
=
\sup_{\substack{f(0)=0\\f(t)=y}}
\left\{
\sum_{\substack{(z,u,w)\in\mathcal P\\z=f(u)}}w
-d\int_0^t|f'(u)|^2\,du
\right\},
\]
where the supremum is over absolutely continuous paths of finite
Dirichlet energy. Our Dirichlet energy is twice the quantity called
entropy by \citet{BT19}.

Write
\[
\mu_D(\nu,d)
=
\lim_{n\to\infty}\frac{\E L_D(\nu,d,(0,n))}{n}.
\]

\begin{theorem}\label{thm: LBCPP}
There exist universal constants $c,\bar c>0$ such that
\[
c\,J_2(\nu)\leq\mu_D(\nu,1)\leq\bar c\,J_2(\nu).
\]
Moreover, if $\nu(y,\infty)=0$, then, for every $t>0$,
\[
\E L_D(\nu,1,(0,t))\geq ctJ_2(\nu)-c'y,
\]
for some $c'>0$.
\end{theorem}

% Entropy-controlled model illustration.
\begin{figure}[tbp]
\centering
% --- requires \usetikzlibrary{arrows.meta,decorations.pathreplacing,calc}
\begin{tikzpicture}[x=1.05cm,y=1.05cm,
    >={Stealth[length=2.2mm]},
    bg/.style={circle, fill=black!28, inner sep=0pt},
    hit/.style={circle, fill=red!75!black, draw=black, line width=0.4pt, inner sep=0pt},
  ]
  % size of a point grows with its weight
  \def\rad#1{\fpeval{0.05+0.038*#1}}

  % light space-time strip [0,t]
  \fill[black!3] (-3,0) rectangle (3,5);
  \draw[black!35, densely dotted] (-3,5) -- (3,5);
  \draw[black!35, densely dotted] (-3,0) -- (3,0);

  % axes
  \draw[->, black!70] (-3.5,-0.3) -- (-3.5,5.5) node[above, black, font=\small] {time $u$};
  \draw[black!70] (-3.58,0) -- (-3.42,0) node[pos=0, left, black, font=\small] {$0$};
  \draw[black!70] (-3.58,5) -- (-3.42,5) node[pos=0, left, black, font=\small] {$t$};
  \draw[->, black!70] (-3.2,-0.6) -- (3.3,-0.6) node[right, black, font=\small] {space $z$};
  \draw[black!70] (0,-0.68) -- (0,-0.52) node[pos=0, below, black, font=\small] {$0$};

  % background Poisson points  (z / u / weight)
  \foreach \z/\u/\w in {0.7/4.37/1.43, -1.54/1.56/6, -2.77/4.01/1.43,
     -1.37/2.24/3.25, 1.64/3.07/1.58, -1.9/3.03/1.05, -2.6/2.57/1.62, 2.34/3.11/1.39,
     -2.73/1.05/1.88, -1.68/1.89/1.01, 1.85/0.88/1.25, 2.33/2.75/2.5, 0.78/3.64/1.09,
     2.08/1.85/1.08, -2.47/1.97/1.04, -1.96/3.99/3.23, 2.68/2.92/2.26, -1.96/2.22/1.74,
     -0.55/0.6/1.97, -1.6/3.31/1.43, 2.49/4.4/1.6, -1.99/1.05/3.6, 0.29/1.0/1.29,
     -0.5/1.28/1.01, -2.63/0.73/2.7, -0.87/2.92/1.16,
     -0.81/2.59/1.05, 2.29/0.86/6, -2.77/3.69/2.23, 1.77/0.22/1.66,
     -1.53/1.08/1.44, -0.9/1.43/1.3, -0.31/4.76/1.24, 1.36/2.88/1.46, 2.37/0.47/1.41,
     -1.39/3.94/1.75, 2.64/1.71/1.45, -1.66/0.39/1.07, 2.33/4.1/1.21, 0.53/3.25/1.15,
     2.58/2.34/1.29, 1.48/3.98/1.41}
     \node[bg, minimum size={2*\rad{\w}cm}] at (\z,\u) {};

  % optimal path: polygonal, straight between the points it collects
  \coordinate (P0) at (0,0);
  \coordinate (P1) at (0.75,0.85);
  \coordinate (P2) at (1.35,1.65);
  \coordinate (P3) at (0.45,2.45);
  \coordinate (P4) at (-0.65,3.15);
  \coordinate (P5) at (-0.45,3.85);
  \coordinate (P6) at (0.1,4.45);
  \coordinate (P7) at (0,5);

  % entropy cost of one segment
  \coordinate (C) at (P2 |- P3);
  \draw[black!55, densely dashed, line width=0.4pt] (P3) -- (C) -- (P2);
  \draw[decorate, decoration={brace, amplitude=3pt}, black!70]
     ($(P3)+(0.12,0.07)$) -- ($(C)+(0,0.07)$) node[midway, above=2pt, font=\scriptsize] {$\Delta z$};
  \draw[decorate, decoration={brace, amplitude=3pt}, black!70]
     ($(C)+(0.07,0)$) -- ($(P2)+(0.07,0.22)$) node[midway, right=3pt, font=\scriptsize] {$\Delta u$};

  \draw[blue!70!black, line width=1.3pt, line join=round]
     (P0) -- (P1) -- (P2) -- (P3) -- (P4) -- (P5) -- (P6) -- (P7);

  % collected points
  \foreach \p/\w/\lab/\pos in {P1/2.5/1/right, P2/5/2/right, P3/1.6/3/below left,
                               P4/4.5/4/left, P5/2/5/left, P6/1.4/6/right}
     \node[hit, minimum size={2*\rad{\w}cm}, label={[font=\scriptsize, red!60!black, inner sep=1pt]\pos:$w_{\lab}$}] at (\p) {};

  % endpoints
  \fill (P0) circle (2.2pt);
  \fill (P7) circle (2.2pt) node[above=2pt, font=\small] {$(0,t)$};

  % value of the path
  \node[anchor=west, font=\small, align=left] at (3.35,2.5)
     {path value\\[3pt]
      $\displaystyle \sum_{i} {\color{red!60!black}w_i}\;-\;d\!\!\sum_{\text{segments}}\!\frac{(\Delta z)^2}{\Delta u}$};
\end{tikzpicture}
\caption{The entropy-controlled Poisson last passage model. Points of a
compound Poisson process carry weights $w_i$ (dot size grows with weight).
A path from $(0,0)$ to $(0,t)$ collects the weights of the points on its
graph and pays the entropy penalty $d\int_0^t f'(u)^2\,du$. An optimal path
is straight between the points it collects, so each segment costs
$d(\Delta z)^2/\Delta u$. The path makes sideways detours to collect
large weights.}
\label{fig:entropy model}
\end{figure}
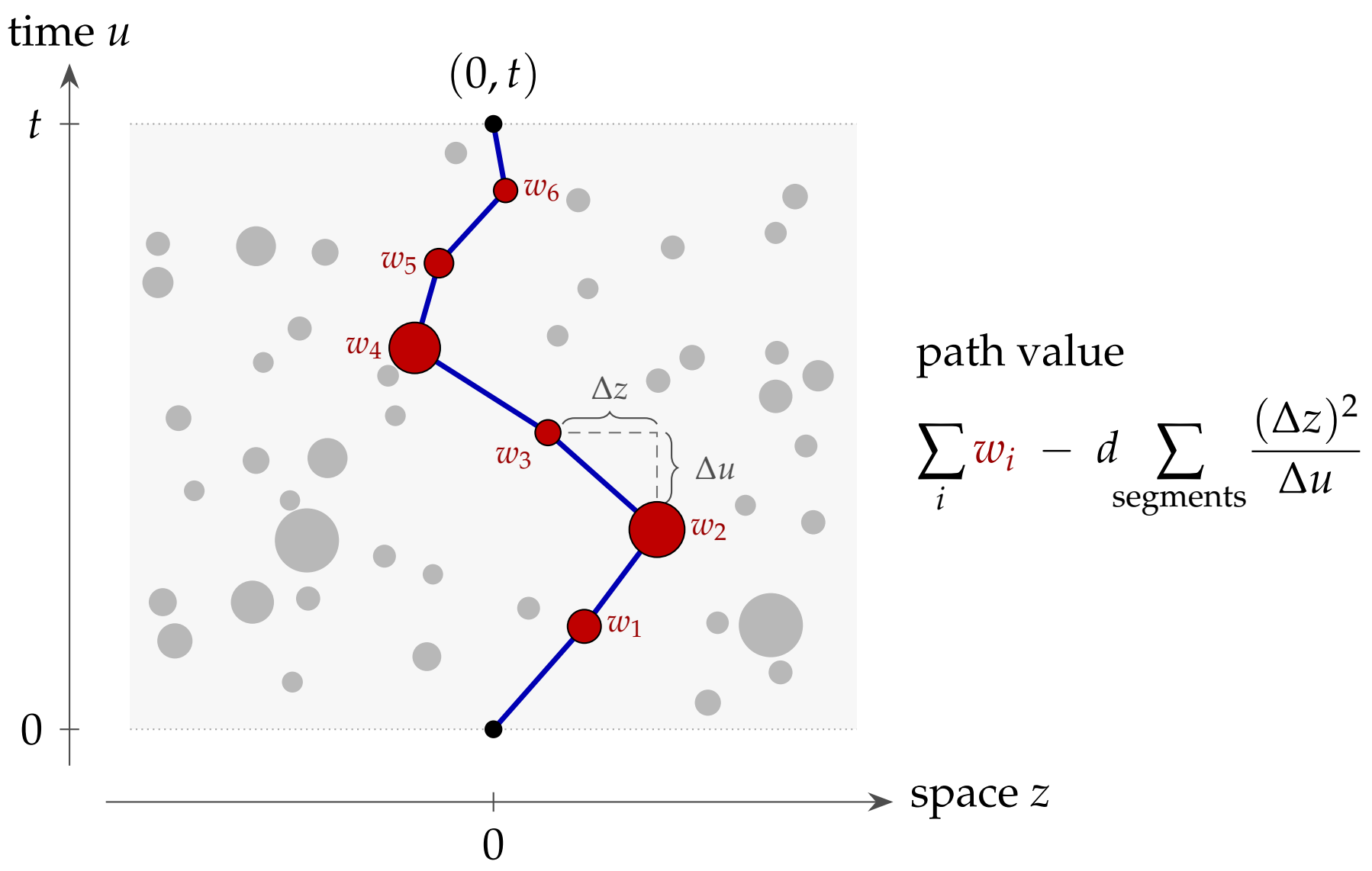

\subsection{Poisson last-passage percolation}

In the ordinary Poisson model, the environment consists of marked
points in the positive quadrant with intensity $dx\,dy\,\nu(dw)$.
Let $L_P(\nu,t)$ be the maximum weight collected by a path from $(0,0)$
to $(t,t)$ whose coordinates are non-decreasing, and define
\[
\mu_P(\nu)=\lim_{t\to\infty}\frac{\E L_P(\nu,t)}{t}.
\]
Unlike in the lattice model, $\nu$ need not be a probability measure.

\begin{theorem}\label{thm: bounds for poisson model}
There exist universal constants $c,\bar c>0$ such that
\[
c\,J_2(\nu)^{3/4}\leq\mu_P(\nu)\leq\bar c\,J_2(\nu)^{3/4}.
\]
\end{theorem}

Thus the lattice and ordinary Poisson models have the same rate,
while the entropy-controlled rate is $J_2(\nu)$ itself.

\subsection{Higher dimensions}

The results and arguments extend to every dimension $q\geq2$.
Let $L_{\Z^q}(\nu,n)$ and $L_{P,q}(\nu,t)$ be the passage values from
the lower to the upper corner of $\{1,\ldots,n\}^q$ and $[0,t]^q$,
respectively. In the entropy-controlled model, paths take values in
$\R^{q-1}$, and the penalty is
\[
d\int_0^t\|f'(u)\|^2\,du.
\]
Denote the corresponding diagonal shape functions by
$\mu_{\Z^q}(\nu)$, $\mu_{P,q}(\nu)$, and $\mu_{D,q}(\nu,d)$. Recall the rate $J_{q}(\nu)$ from \eqref{eq: gen dim rate}.

\begin{theorem}\label{thm: gen dim bounds}
For every $q\geq2$, there exist constants $c_q,\bar c_q>0$, depending
only on $q$, such that
\[
\begin{aligned}
c_qJ_q(\nu)^{(q+1)/(2q)}
&\leq\mu_{\Z^q}(\nu)
\leq\bar c_qJ_q(\nu)^{(q+1)/(2q)},\\
c_qJ_q(\nu)^{(q+1)/(2q)}
&\leq\mu_{P,q}(\nu)
\leq\bar c_qJ_q(\nu)^{(q+1)/(2q)},\\
c_qJ_q(\nu)
&\leq\mu_{D,q}(\nu,1)
\leq\bar c_qJ_q(\nu).
\end{aligned}
\]
For the lattice statement, $\nu$ is a probability measure on $[0,\infty)$.
\end{theorem}

In particular, the same tail integral characterizes finiteness in all
three models.

\subsection{Critical heavy tails}

The finite-size estimates also determine the growth rate when the
shape function is infinite. \citet{hambly-martin} established scaling
limits for two-dimensional last-passage percolation with regularly
varying tails of exponent $\alpha<2$. At the critical exponent
$\alpha=2$, contributions from many weight scales accumulate.

\citet{ganguly2024last} already obtained the correct logarithmic
exponent in their lower bound. For critical inverse-square tails,
they proved, with high probability,
\[
L_{\Z^2}(\nu,n)\geq
c\,\frac{n(\log n)^{3/4}}{\log\log n}.
\]
The  factor $\log\log n$ arises from
the separation between the scales used to obtain the necessary
concentration estimates. Their higher-dimensional lower bound
likewise contains the exponent $(q+1)/(2q)$, again with a
$\log\log n$ loss.

Our finite-size bounds establish that these logarithmic exponents
are sharp.

\begin{theorem}\label{thm: critical heavy tailed LPP bound}
Suppose $q\geq2$ and
\[
\nu(x,\infty)\sim c_0x^{-q}
\qquad\text{as }x\to\infty,
\]
for some $c_0>0$. There exist constants $c,\bar c>0$, depending on $q$
and $\nu$, such that, for every $n\geq2$,
\[
c\,n(\log n)^{(q+1)/(2q)}
\leq\E L_{\Z^q}(\nu,n)
\leq\bar c\,n(\log n)^{(q+1)/(2q)}.
\]
\end{theorem}

In dimension two, the mean passage value therefore has order
$n(\log n)^{3/4}$. The theorem removes the iterated-logarithmic loss
from the lower bound for the mean and supplies a matching upper bound,
determining the growth up to constant factors.

\section{Scaling and the rate integral}

\begin{lemma}[Scaling and shear]\label{lem:entropy scaling}
Let $q\geq2$ and $M_h\nu(A)=\nu(A/h)$. For $a,h,d,t>0$, put
$s=a^{2/(q+1)}(h/d)^{(q-1)/(q+1)}$. Then, for
$\mathbf y\in\R^{q-1}$,
\[
L_{D,q}(aM_h\nu,d,(\mathbf y,t)) \stackrel{d}{=}hL_{D,q}(\nu,1,(\mathbf0,st)) -\frac{d\|\mathbf y\|^2}{t},\qquad \mu_{D,q}(aM_h\nu,d)=hs\,\mu_{D,q}(\nu,1).
\]
\end{lemma}
\begin{proof}
The shear $(z,u)\mapsto(z-u\mathbf y/t,u)$ preserves the Poisson
{intensity and reduces the Dirichlet energy by $\|\mathbf y\|^2/t$,
and hence the penalty by $d\|\mathbf y\|^2/t$.}
For the zero endpoint, use $(z,u,w)\mapsto(\xi z,\tau u,hw)$,
where $\xi=(h/(ad))^{1/(q+1)}$ and $\tau=1/s$.
Since $\xi^{q-1}\tau=a^{-1}$ and $d\xi^2/\tau=h$, the transformed
intensity is $aM_h\nu$ and every path value is multiplied by $h$.
The shape relation follows by taking expectations, dividing by $t$,
and letting $t\to\infty$.
\end{proof}

\begin{lemma}[Sum-integral comparison]\label{lem: rate comparison}
Let $q\geq2$, $r>1$, and $a_i=\nu[r^i,r^{i+1})$. There exist
constants $c_{q,r},\bar c_{q,r}>0$ such that
\[
c_{q,r}J_q(\nu)
\leq\sum_{i\in\Z}a_i^{2/(q+1)}r^{2qi/(q+1)}
\leq\bar c_{q,r}J_q(\nu).
\]
\end{lemma}
\begin{proof}
Write $f(x)=\nu(x,\infty)$, $p=2/(q+1)$, and
$\alpha=2q/(q+1)$. For $r^{i-1}\leq x<r^i$, we have
$f(x)\geq a_i$, and hence
\[
\int_{r^{i-1}}^{r^i}f(x)^p x^{\alpha-1}\,dx
\geq c_{q,r}a_i^p r^{\alpha i}.
\]
This gives the second inequality. Conversely, for
$r^i\leq x<r^{i+1}$, we have $f(x)\leq\sum_{j\geq i}a_j$.
Since $p\leq1$,
\[
J_q(\nu) \leq\bar c_{q,r}\sum_i r^{\alpha i} \left(\sum_{j\geq i}a_j\right)^p \leq\bar c_{q,r}\sum_j a_j^p\sum_{i\leq j}r^{\alpha i} \leq\bar c_{q,r}\sum_j a_j^p r^{\alpha j}.\qedhere
\]
\end{proof}

\section{Proofs for the entropy controlled LPP model}

We first derive a lower bound for $\mathbb{E} L_D(a \delta_b, 1, (0, t))$ in the entropy-controlled Poisson last passage model. For $\nu = a \delta_b$, intensity is $a$, weights are $b$. We use the bound $\mathbb{E} L_D(\delta_1, 1, (0, t)) \geq c \min(t, t^{3/2})$.
\begin{proposition}\label{prop: basic lb}
For the entropy-controlled Poisson last passage model with measure $\nu = a \delta_b$, where $\delta_b$ is the Dirac delta at $b > 0$, we have:
\[
\mathbb{E} L_D(a \delta_b, 1, (0, t)) \geq c a^{2/3} b^{4/3} \min(t, a^{1/3} b^{1/6} t^{3/2}) \ge ca^{2/3}b^{4/3}t - cb.
\]
\end{proposition}
\Needspace{9\baselineskip}
% \begin{revision}
\begin{remark}
The crossover time is of order $d^{1/3}a^{-2/3}b^{-1/3}$. More precisely,
for $t\geq d^{1/3}a^{-2/3}b^{-1/3}$ and $x\in\R$,
\[
c\,t\mu_D(a\delta_b,d)
\leq \E L_D(a\delta_b,d,(x,t))+\frac{dx^2}{t}
\leq t\mu_D(a\delta_b,d),
\]
with a universal constant $c>0$. Thus the expected value is comparable
to its linear growth rate after the deterministic shear cost is added
back. For $d=1$, the crossover time is $a^{-2/3}b^{-1/3}$.
\end{remark}
% \end{revision}
\begin{proof}[Proof of Proposition \ref{prop: basic lb}]
For $s\leq1$, a point in
$[-\sqrt{s}/4,\sqrt{s}/4]\times[s/3,2s/3]$ gives a passage value
bounded below by a positive constant. The probability of such a point
is at least $cs^{3/2}$. For $s\geq1$, concatenation of unit intervals
gives a lower bound $cs$. Thus, with $s=a^{2/3}b^{1/3}t$,
Lemma \ref{lem:entropy scaling} gives
\[
\E L_D(a\delta_b,1,(0,t))
=b\E L_D(\delta_1,1,(0,s))
\geq cb\min(s,s^{3/2}).
\]
This is the first bound, and $\min(s,s^{3/2})\geq s-1$
gives the second.
\end{proof}
\begin{corollary}\label{c:EL}
% \begin{revision}
Let $p_i=(x_i,t_i)\in\R^2$, with $t_0<\cdots<t_n$, and let
\[
d(p)=\sum_{i=1}^n\frac{(x_i-x_{i-1})^2}{t_i-t_{i-1}}
\]
be the minimum Dirichlet energy of a path passing through all the $p_i$.
% \end{revision} 
Let $L_D(a\delta_b,1,p)$ be the value of the best last passage path through $p_i$. Then
$$
\mathbb{E} L_D(a\delta_b,1,p) \geq c a^{2/3} b^{4/3} (t_n-t_0) - c b n - d(p).
$$
\end{corollary}
\begin{proof}
By Lemma \ref{lem:entropy scaling} and the bound above, we have
\[
\mathbb{E} L_D(a\delta_b,1,p)+d(p) = \sum_{i=1}^n \mathbb{E} L_D(a\delta_b,1,(0,t_i-t_{i-1})) \geq c a^{2/3} b^{4/3} (t_n-t_0) - c b n.\qedhere
\]
\end{proof}

We now prove Theorem \ref{thm: LBCPP}. The approach for the lower bound constructs a last passage path for a discrete measure and extends the result to general measures using stochastic domination and tail integrals. The proof of the upper bound is based on first considering the largest weight and obtaining an upper bound on the exponential of passage time over every possible path. Then conditioning on higher weights a similar bound is obtained on the lower weights. 

\Needspace{7\baselineskip}
\begin{proof}[Proof of Theorem \ref{thm: LBCPP}]
% Step 1: recursive path construction for discrete measure
\hspace{0.1mm}
\subsubsection*{{Proof of the lower bound}}
We first consider integers $m\le n$ and discrete measure \(\nu = \sum_{i=m}^n a_i \delta_{b_i}\), where \(a_i ,b_i \ge 0 \). We construct a last passage path recursively to bound \(\mathbb{E} L_D(\nu, 1, (0, t))\).
\subsubsection*{ Recursive path construction}
Initialize with \(P_{n+1} = \{(0, 0), (0, t)\}\), \(N_{n+1} = 2\). For \(i = n, n-1, \ldots, m\), let $b_i'\ge 0$ to be specified later. For each $i$, consider independent Poisson processes $R_i$ with intensity $a_i$.
At step \(i\):
\begin{enumerate}[label=(\alph*)]
    \item Consider a path $f$ that passes through all points of $P_{i+1}$ and maximizes
    $$
    b_i'\cdot\#\{(z, u)\in R_i: z = f(u)\} - \int_0^t [f'(u)]^2 \, du
    $$
   
    \item Let \(F_i\) be the set of $R_i$ points on $f$, and set \(N_i = |F_i|\).
    \item Update \(P_i = P_{i+1} \cup F_i\), and let \(\mathcal{F}_i\) be the \(\sigma\)-field generated by $R_i,\ldots, R_n$.
\end{enumerate}

The maximizing paths and their point counts are well defined. Conditional
on the parent polygon, Corollary~\ref{cor:polygon} below gives finite
expectation for the exponential sum over all refinements. Almost surely,
only finitely many refinements have nonnegative score relative to the
parent, and the empty refinement has score zero. A maximizer therefore
exists; fix a measurable tie-breaking rule. For a maximizing refinement
with reward $b>0$, point count $N$, and excess energy $\Delta E$, we have
$0\leq\Delta E\leq bN$. Evaluating that refinement with reward $2b$ gives
\[
bN\leq 2bN-\Delta E
\leq\theta^{-1}\log W_{\theta}^{(2b)}(\sigma),\qquad \theta=b^{-1},
\]
where $W_{\theta}^{(2b)}(\sigma)$ is the insertion sum of
Corollary~\ref{cor:polygon} with point weight $2b$.
Conditional Jensen and that corollary bound the right-hand side in
expectation uniformly in the parent vertices. Induction gives integrable
point counts and energies at every step, justifying the expectations
and telescoping below. Zero-intensity layers may be omitted.

Using the lower bound from Corollary \ref{c:EL} we estimate
\[
\mathbb{E}[b_i' N_i - d(P_i) | \mathcal{F}_{i+1 }] \geq c t a_i^{2/3}b_i'^{4/3} - c b_i' \#P_{i+1} - d(P_{i+1}),
\]
where
$$d(P_i) = \sum_{k=1}^{|P_i|-1} \frac{(z_{k+1} - z_k)^2}{u_{k+1} - u_k}$$ is the {Dirichlet energy}, \(\#P_{i+1} = \sum_{j=i+1}^{n+1} N_j\).
We take expectations and sum over $i$ to get
$$
\sum_{i=m}^n b_i' \E N_i + \sum_{i=m}^n \E[d(P_{i+1})-d(P_i)] \ge ct\sum_{i=m}^n a_i^{2/3} {b_i'}^{4/3} -c\sum_{m\le i< j \le n+1} b_i'\E N_{j}
$$
Now let $s_i=b_i'+c\sum_{j=m}^{i-1} b_j'$, $s_*=\sum_{j=m}^{n}b_j'$. Note also that the second sum is telescoping, $d(P_{n+1})=0$ and $N_{n+1}=2$. So we get
$$
\E\l[\sum s_iN_i - d(P_m)\r] \ge ct\sum_{i=m}^na_i^{2/3} {b_i'}^{4/3}-2cs_*.
$$
We now set $r=2c+1$, and $b_i=r^i, b_i'=b_i/2$. Then
\[
s_i\leq b_i,\qquad 2cs_*=\frac{r^{n+1}-r^m}{2},
\]
and we get
$$
\E L_D(\nu,1,(0,t)) \ge 2^{-4/3}ct\sum_{i=m}^n a_i^{2/3}b_i^{4/3}-\frac{1}{2}r^{n+1}.
$$
\subsubsection*{ General $\nu$}
Now for a general measure $\nu$ supported on $[0,y]$, with $y>0$,
choose $n$ so that $r^n\leq y<r^{n+1}$ and put
\[
a_i=\nu[r^i,r^{i+1}),\qquad
\nu_r=\sum_{i=m}^n a_i\delta_{r^i}.
\]
Then $\nu$ stochastically dominates $\nu_r$. By the preceding bound
and Lemma \ref{lem: rate comparison},
\[
\E L_D(\nu,1,(0,t))
\geq c'tJ_2\bigl(\nu|_{[r^m,r^{n+1})}\bigr)-c''y.
\]
Letting $m\to-\infty$ gives the second claimed bound. Dividing by $t$
and letting $t\to\infty$ gives the shape function lower bound for
$\nu$ with bounded support. For $\nu$ with unbounded support, apply this bound to its restriction
to $[0,x]$ and let $x\to\infty$, using monotonicity and the monotone
convergence theorem.

We now prove the upper bound. We use Lemma \ref{lem:insert} which is proved later.

\subsubsection*{Proof of upper bound}
Fix $(x,0)$ and $(y,t)$, and let $\Pi_a$ be a Poisson process of intensity $a\dd z\dd s$ on $\R\times(0,t)$ with each point of weight $w$. A fixed path $\pi$ is a sequence $\pi=(p_1,\ldots,p_m)$ of finitely many points in strictly increasing time order. {Set $p_i=(z_i,s_i)$ and add endpoints $(z_0,s_0)=(x,0)$, $(z_{m+1},s_{m+1})=(y,t)$.} Define
\begin{equation}\label{eq:score}
{E_{x,y;t}(\pi)=\sum_{i=0}^m\frac{(z_{i+1}-z_i)^2}{s_{i+1}-s_i}},
\qquad
\Delta E(\pi)=E_{x,y;t}(\pi)-\frac{(y-x)^2}{t}.
\end{equation}

For $\theta>0$, set
\begin{equation}\label{eq:W}
W_{a,w,\theta}(x,y;t)
=\sum_{m=0}^\infty\ 
\sum_{\substack{\pi=(p_1,\ldots,p_m)\\0<s_1<\cdots<s_m<t}}
\exp\{\theta mw-\theta \Delta E(\pi)\}.
\end{equation}
For $m=0$ the inner sum is the single term $1$. The above sum is over every possible path $\pi$. The expectation in the following lemma is only over $\Pi_a$.

\begin{lemma}\label{lem:insert}
For every $x,y\in\R$ and $a,w,\theta,t>0$,
\begin{equation}\label{eq:insertbound}
\E W_{a,w,\theta}(x,y;t)
\le \exp\{\bar c_0a^{2/3}\theta ^{-1/3}e^{2\theta w/3}t\},
\qquad\text{where } \bar c_0=\frac32\pi^{1/3}.
\end{equation}
\end{lemma}

Let $\sigma$ be a polygonal path of time duration $t$ and $\tau$ be a refinement of $\sigma$ obtained by adding points while retaining vertices of $\sigma$. Suppose the new added points are from a Poisson point process of intensity $a$ and weight $w$ and let $N_{\text{new}}(\tau)$ denote the new Poisson points added. Define
\[
W_{\theta}(\sigma)=\sum_{\tau\text{ refines }\sigma}
e^{\theta\{wN_{\rm new}(\tau)-[E(\tau)-E(\sigma)]\}}
\]
Lemma \ref{lem:insert} immediately gives the following corollary.
\begin{corollary}\label{cor:polygon}
\begin{equation}\label{eq:polygonbound}
\E W_{\theta}(\sigma)\le
\exp\{\bar c_0a^{2/3}\theta ^{-1/3}e^{2\theta w/3}t\}.
\end{equation}
The bound is uniform in the number and locations of the vertices of $\sigma$.
\end{corollary}
% \begin{revision}
\begin{proof}
The insertion sum factorizes over the segments of $\sigma$.
The corresponding Poisson processes are independent because their time
intervals are disjoint. Multiplying the bounds in Lemma~\ref{lem:insert}
and adding the segment durations proves the claim.
\end{proof}
% \end{revision}

    We first prove the upper bound when $\nu=\sum_{j=1}^na_j\delta_{w_j}$ where $a_j>0$ and $w_1>w_2>\dots>w_n>0$. Let $\theta_j= w_j^{-1}$. Let $\mathcal{S}_i$ contain all paths (from $(0,0)$ to $(0,t)$) of finitely many Poisson points and with weights from $w_j$ for $j\leq i$. For a path $\pi$, let $L_D(\pi)$ denote the passage time in the entropy controlled model, which is the sum of the weights on the path $\pi$ with its entropy penalty reduced. Define 
    \[
    Z_i=\sum\limits_{\sigma\in \mathcal{S}_i}e^{\theta_i L_D(\sigma)}, \quad F_i=\frac{1}{\theta_i}\log Z_i.
    \]
    We have included the empty path in $\mathcal{S}_i$ so we have $Z_i\geq 1$ and put $F_0=0$.
    Deleting the points of weight $w_i$ from any path in $\mathcal{S}_i$ gives a path in $\mathcal{S}_{i-1}.$ Recall the notation of $W_{\theta}(\sigma)$ from Corollary \ref{cor:polygon}. It can be checked easily that
    \[
    Z_i=\sum\limits_{\sigma\in \mathcal{S}_{i-1}}e^{\theta_i L_D(\sigma)}W_{\theta_i}(\sigma).
    \]
    Let $\mathcal{F}_i$ be the sigma algebra generated by the Poisson points of weight $w_j$ for $j\leq i$. Applying Corollary \ref{cor:polygon} we get that
    \[
    \E[Z_i|\mathcal{F}_{i-1}]\leq e^{\theta_i A_it}\sum\limits_{\sigma\in \mathcal{S}_{i-1}}e^{\theta_i L_D(\sigma)},
    \]
    where $A_i=\bar c_*a_i^{2/3}w_i^{4/3}$ for some constant $\bar c_*.$ Conditional Jensen inequality gives us that
    \[
    \E[F_i|\mathcal{F}_{i-1}]\leq A_it+\frac{1}{\theta_i}\log\sum\limits_{\sigma\in\mathcal{S}_{i-1}}e^{\theta_i L_D(\sigma)}.
    \]
{The same calculation as in the proof of Lemma \ref{lem:insert},
with $ae^{\theta w}$ replaced by $\sum_{j=1}^{i}a_j e^{\theta w_j}$, shows that}
$\mathbb{E}Z_i<\infty$. In particular, $\log Z_i$ is integrable,
and conditional Jensen's inequality is applicable.
    {For $i\geq2$, use $\sum_j x_j^{p_i}\leq(\sum_j x_j)^{p_i}$
with $p_i=\theta_i/\theta_{i-1}\geq1$ to obtain}
\[
{\sum_{\sigma\in\mathcal S_{i-1}}e^{\theta_iL_D(\sigma)}\leq Z_{i-1}^{p_i}.}
\]
Hence we have \[
    \E[F_i|\mathcal{F}_{i-1}]\leq A_it+F_{i-1}.
    \]
    Note that this holds for $i=1$ with $F_0=0$. Taking expectations and iterating this process we get 
    \[
    \E F_n\leq t\sum\limits_{i=1}^nA_i.
    \]
    {For every $\sigma\in\mathcal S_n$, $Z_n\geq\exp(\theta_n L_D(\sigma))$.
Taking the supremum gives $F_n\geq L_D(\nu,1,(0,t))$.} Taking expectations we get 
    \[
    \E L_D\l(\nu,t\r)\leq  \bar c_*t\sum\limits_{j=1}^n a_j^{2/3}w_j^{4/3}.
    \]

    For general $\nu$, take $a_k=\nu[2^{k},2^{k+1})$ for $k\in\Z$. Restrict the weights to $[2^{-n},2^n)$ for $n$  large and fixed. Rounding the weights in $[2^{k},2^{k+1})$ upward to $2^{k+1}$, we construct measures $\nu_n:=\sum_{i=-n}^{n-1}a_i\delta_{2^{i+1}}$ and $\Tilde{\nu}_n:=\nu|_{[2^{-n},2^n)}$.  

    By the above proof and Lemma \ref{lem: rate comparison},
    \[
\E L_D\l(\Tilde{\nu}_n,t\r) \leq\E L_D\l(\nu_n,t\r) \leq\bar c t\sum_{j=-n}^{n-1}a_j^{2/3}2^{4j/3} \leq\bar c tJ_2(\nu).
\]
    Letting $n\to\infty$, monotone convergence gives
    $\E L_D(\nu,t)\leq\bar c tJ_2(\nu)$.
    Taking $t\to\infty$ completes the proof of the upp  er bound for general $\nu$.
\end{proof}

\begin{proof}[Proof of Lemma \ref{lem:insert}]
We write $W_{a,w,\theta}(x,y;t)$ as just $W$ in the proof.
The factorial-moment formula for Poisson processes, together with Tonelli's theorem, gives the exact identity
\begin{equation}\label{eq:campbell}
\E W=\sum_{m=0}^\infty a^m e^{\theta mw}
\int_{0<s_1<\cdots<s_m<t}\int_{\R^m}
e^{-\theta \Delta E}\dd\mathbf z\,\dd\mathbf s.
\end{equation}
As the tuple times are already restricted to be increasing, there is no extra $1/m!$ in this formula. For the Poisson formula, see \citet[Theorem 4.4 and Corollary 4.10]{LP18}.

To bound the integral, we rewrite $\Delta E$.
Let $\ell(s)=x+(y-x)s/t$ and put $u_i=z_i-\ell(s_i)$, with $u_0=u_{m+1}=0$. Write $s_0=0$, $s_{m+1}=t$ and $\Delta_i=s_{i+1}-s_i$. It is easy to check that
\begin{equation}\label{eq:energy}
\Delta E=\sum_{i=0}^m\frac{(u_{i+1}-u_i)^2}{\Delta_i}.
\end{equation}
The change of spatial variables has Jacobian one. As a result,
\[
\int_{\R^m}e^{-\theta\Delta E}\,d\mathbf z
=\int_{\R^m}e^{-\theta\Delta E}\,d\mathbf u.
\]
Gaussian integration gives us,
\[
\int\limits_{\R^m}e^{-\theta \Delta E }d\mathbf{u}=\l(\frac{\pi}{\theta }\r)^{m/2}\sqrt{\frac{\Delta_0\cdots\Delta_m}{t}}
\]
Since $\sum_i\Delta_i=t$, we have that \[
\l(\Delta_0\cdots\Delta_m\r)^{1/(m+1)}\leq \frac{t}{m+1}.
\]
{Using this, the time-simplex volume $t^m/m!$, and
$m!\leq(m+1)^{m+1}$, we get}
\[
\E W\leq \sum\limits_{m=0}^{\infty}\frac{v^m}{\l(m!\r)^{3/2}}\leq e^{(3/2)v^{2/3}},
\]
where
\[
v=ae^{\theta w}\sqrt{\frac{\pi}{\theta}}\,t^{3/2}.
\]
% \begin{revision}
% For the last inequality, put $b_m=v^{2m/3}/m!$ and use
% \[
% \sum_{m\geq0}\frac{v^m}{(m!)^{3/2}}
% =\sum_{m\geq0}b_m^{3/2}
% \leq\left(\sum_{m\geq0}b_m\right)^{3/2}
% =e^{(3/2)v^{2/3}}.
% \]
% \end{revision}
This gives the desired bound.
\end{proof}

\section{Proofs of Poisson LPP model}
\begin{proof}[Proof of Theorem \ref{thm: bounds for poisson model}]
The upper bound in Theorem \ref{thm: bounds for poisson model} follows from the upper bound in Theorem \ref{thm: LBCPP}. Rotate points $(x,y)$ in the Poisson model to the points
\[
{(u,v)=\left(\frac{x-y}{\sqrt2},\frac{x+y}{\sqrt2}\right)}
\]
in the entropy model, for some parameter $d$. The transformation is area preserving and $(0,0),(n,n)$ are mapped to $(0,0),(0,\sqrt{2}n)$ and {the rotated path has Dirichlet energy at most $\sqrt{2}n$,
and hence penalty at most $d\sqrt{2}n$.} This coupling gives,
\[
L_P(\nu,n)\leq L_D(\nu,d,(0,\sqrt{2}n))+\sqrt{2}dn,
\]
for any parameter $d>0$.

After dividing by $n$, taking $n\rightarrow\infty$, and applying Theorem \ref{thm: LBCPP} and Lemma \ref{lem:entropy scaling}, optimization over $d$ gives,
\[\mu_P(\nu) \leq \bar c \l(\int_0^\infty \nu(x, \infty)^{2/3} x^{1/3} \, dx\r)^{3/4}.\]

In order to prove the lower bound in Theorem \ref{thm: bounds for poisson model}, we first prove the following theorem.

\begin{theorem}\label{thm: lower bound for poisson model}
Let $\nu$ be a measure on $\mathbb{R}_+$. Let $b_i=2^i$ and $\nu[2^{i},2^{i+1})=a_i$. Fix $k$ and $\nu_{\ell}=\sum_{n=1}^ka_n\delta_{b_n}$. There exists $c_1>0$(independent of $k$) such that for any choice of $\varepsilon_n\in[0,1],$ we have
\[
\mu_P(\nu_{\ell})\geq c_1\prod\limits_{n=1}^k(1-\varepsilon_{n}^2)\sum_{n=1}^{k}b_{n}\sqrt{a_{n}\varepsilon_{n}/2}.
\]
\end{theorem}

The proof of Theorem \ref{thm: lower bound for poisson model} is based on a recursive path construction argument described as follows. We pick $\{\varepsilon_n\}_{n=1}^{k}$, $\varepsilon_n\in (0,1)$ (to be decided later). We construct a path of line segments recursively. We are given a path after stage $n$ (chronologically step indexed by $k$ is performed first and then $k-1$ and so on), which goes through vertices of weight at least $b_{n}$. This set of vertices is denoted $E_n$, and is ordered naturally so that one obtains an up-right path. Initially, i.e. after stage $k+1$, the set of vertices is $E_{k+1} = \{ (0,0), (t,t)\}$.

At stage $n$, we run Poisson last passage on $\nu = a_{n} \delta_{b_{n}}$, between every consecutive pair of points in $E_{n+1}$, with an additional constraint; if the slope between two consecutive points in $E_{n+1}$ is $m$, then the additional constraint is that the slope of consecutive line segments must be between $m(1-\varepsilon_{n})$ and $m/(1-\varepsilon_{n})$. If the points picked up by this last passage are denoted by $P_{n}$, then $E_{n} = P_{n} \cup E_{n+1}$, and then we proceed to step $n-1$.
Figure~\ref{fig:multiscale} shows this construction on a simulated example.

\begin{figure}[tbp]
\centering
\includegraphics[width=0.92\textwidth]{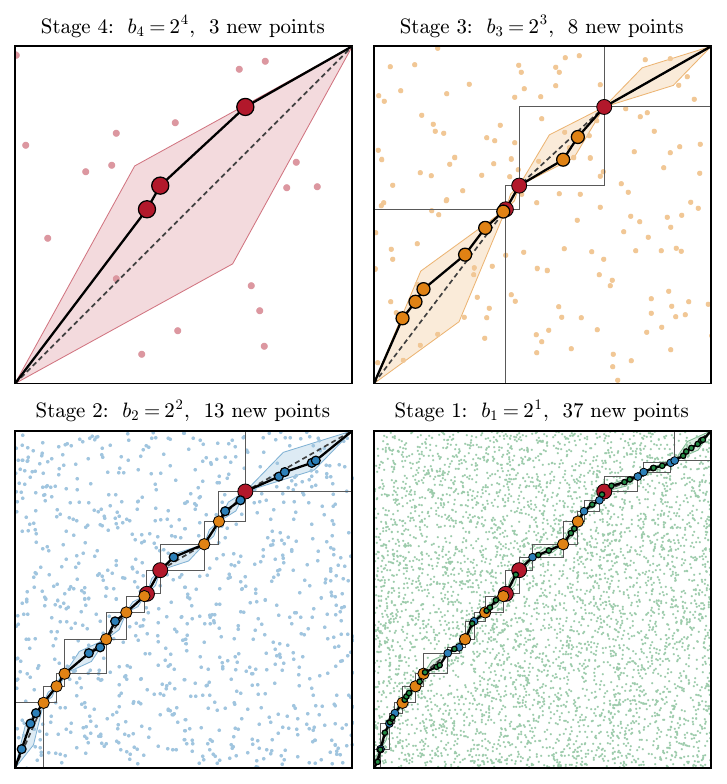}
\caption{The multiscale path construction in the proof of
Theorem~\ref{thm: lower bound for poisson model}, on a simulated Poisson
process with four weight scales $b_n=2^n$ (heavier points are rarer).
At stage $n$, last passage percolation on the points of weight $b_n$ is run
inside each rectangle $R_{n+1,i}$ spanned by consecutive points of $E_{n+1}$
(boxes), restricted to the shaded parallelogram of slopes
$m(1-\varepsilon_n)^{\pm1}$. The dashed line is the path after stage $n+1$;
the solid line is the refined path. The few heaviest points fix the
macroscopic course of the path, and each lighter scale is threaded along it
on a finer scale, each contributing order
$t\,b_n\sqrt{a_n\varepsilon_n}$ to the passage time.
Here $\varepsilon_n=0.45$ for visibility.}
\label{fig:multiscale}
\end{figure}

We first compute a lower bound on {$\E L_P(\sum_{i=1}^{k} a_i \delta_{b_i},t)$}. 
The points in $E_i$ are totally ordered by the coordinates; let $p_s=(x_{i,s},y_{i,s})$ be the $s+1$th point in $E_i$. The points $(x_{i,0},y_{i,0})$ and $(x_{i,|E_i|-1},y_{i,|E_i|-1})$ also belong to $E_{i-1}$. 
    Define the rectangle $R_{i,s}=[x_{i,s-1},x_{i,s}]\times [y_{i,s-1},y_{i,s}]$ with area $A_{i,s}$ for $s\in[|E_i|-1]$.

    Since $E_{k+1}= \{ (0,0), (t,t)\}$, at stage $k$, between these two points in $E_{k+1}$, we perform Poisson LPP on $\nu=a_k\delta_{b_k}$ with additional constraint that {slopes of the new consecutive segments lie between $1-\varepsilon_k$
and $1/(1-\varepsilon_k)$.}

{Write $\lambda(s)=\E L_P(\delta_1,s)$.}
We make use of the following lemmas, the proofs of which are postponed.

% \begin{revision}
\begin{lemma}\label{lem: delta_1 lb}
There exists a universal constant $c_1>0$ such that
\[
\E L_P(a\delta_b,1)\geq b(c_1\sqrt a-1)
\]
for all $a,b>0$.
\end{lemma}
% \end{revision}

\begin{lemma} \label{smallerrectangle}
    For any subrectangle $R_{n-1,j} \subset R_{n,i}$, define relative length and width by
\[
l_j=\frac{x_{n-1,j}-x_{n-1,j-1}}{x_{n,i}-x_{n,i-1}},
\qquad
w_j=\frac{y_{n-1,j}-y_{n-1,j-1}}{y_{n,i}-y_{n,i-1}}.
\]
Then

    \begin{equation*}
        A_{n-1,j} \geq (l_j+w_j)^2\frac{(1-\varepsilon_{n-1})^2}{((1-\varepsilon_{n-1})^2+1)^2}A_{n,i}.
    \end{equation*}
\end{lemma}

\begin{lemma} \label{smallerrectangles}
    For any rectangle $R_{n,i}$, we have

    \begin{equation*}
        \sum_{R_{n-1,j} \subset R_{n,i}} \sqrt{A_{n-1,j}} \geq \left(1-{\varepsilon_{n-1}^2}\right)\sqrt{A_{n,i}}
    \end{equation*}
\end{lemma}

\begin{lemma}\label{parallelogram}
Define $l_{\varepsilon}(A\nu)$ to be the last passage time in a rectangle of area $A$ with intensity measure $\nu$, with the additional constraint that all lines have slope between $m(1-\varepsilon)$ and $m/(1-\varepsilon)$, where $m$ is the slope of the rectangle. Then
\begin{equation}\label{eq: rectangle passage time relation}
    l_{\varepsilon}(A\nu) \stackrel{d}{=} L(\nu,\sqrt{\alpha(\varepsilon)A}),
\end{equation}

where,

\begin{equation*}
\alpha(\varepsilon)=\frac{\varepsilon}{2-\varepsilon}\geq \frac{\varepsilon}{2}, \qquad \varepsilon \in [0,1]
\end{equation*}

is the area of the parallelogram with opposite vertices $(0,0)$ and $(1,1)$ and the slopes of the sides are $(1-\varepsilon)^{\pm1}$.
\end{lemma}

\begin{proof}[Proof of Theorem \ref{thm: lower bound for poisson model}]

Let $\mathcal F_n$ be the $\sigma$-field generated by all points of weight at least $b_n$. For $i\in[|E_n|-1]$, let $X_{n-1,i}$ be the number of points $P_{n-1}\cap R_{n,i}$.
We have 
$$
\E[X_{n-1,i}|\mathcal F_n]=\lambda\l(\sqrt{a_{n-1}A_{n,i}\alpha(\varepsilon_{n-1})}\r),
$$

giving the bound
% \begin{revision}
\begin{align}\label{eq: lb for number of points}
\E[X_{n-1,i}\mid\mathcal F_n]
&\geq\lambda\left(\sqrt{a_{n-1}A_{n,i}\varepsilon_{n-1}/2}\right)\notag\\
&\geq c_1\sqrt{a_{n-1}A_{n,i}\varepsilon_{n-1}/2}-1
\end{align}
% \end{revision}
by Lemma \ref{parallelogram} and Lemma \ref{lem: delta_1 lb}. Indeed, since conditional on $\mathcal{F}_n,$ between two consecutive points of $E_{n}$, the rectangle has area $A_{n,i}$. The constraint on the slope (should be between $m_i(1-\varepsilon_{n-1})$ and $m_i/(1-\varepsilon_{n-1})$) for points picked in the next stage, that is the points $P_{n-1} \cap R_{n,i}$, is same as performing last passage percolation on points of weight $b_{n-1}$ which occur at rate $a_{n-1}$ but on rectangle of area $A_{n,i}\alpha(\varepsilon_{n-1})$.

% \begin{revision}
Let $L_n=b_n|P_n|$ be the contribution of level $n$ to the
constructed path, and put $W=\sum_{n=1}^kL_n$.
By \eqref{eq: lb for number of points} and Lemma~\ref{smallerrectangles},
\[
\E[L_{n-1}\mid\mathcal F_n]
\geq c_1t\,b_{n-1}\sqrt{a_{n-1}\varepsilon_{n-1}/2}
\prod_{j=n}^k(1-\varepsilon_j^2)-b_{n-1}(|E_n|-1).
\]
Since $|E_{n+1}|-1=1+\sum_{m>n}|P_m|$ and
$\sum_{n<m}b_n\leq b_m$, the total error satisfies
\[
\sum_{n=1}^k b_n(|E_{n+1}|-1)
=\sum_{n=1}^k b_n+\sum_{m=1}^k|P_m|\sum_{n<m}b_n
\leq2b_k+W.
\]
Taking expectations and summing the conditional bounds therefore gives
\[
2\E W\geq c_1t\prod_{n=1}^k(1-\varepsilon_n^2)
\sum_{n=1}^k b_n\sqrt{a_n\varepsilon_n/2}-2b_k.
\]
Since $L_P(\nu_\ell,t)\geq W$, division by $t$ and passage to the limit
prove the theorem, after changing $c_1$.
If any $\varepsilon_n=1$, the claimed lower bound is zero; levels with
$\varepsilon_n=0$ are simply left unused.
% \end{revision}
\end{proof}

\Needspace{11\baselineskip}
\begin{samepage}
Taking
\[
\epsilon_n=\frac{b_n^{2/3}a_n^{1/3}}
 {\sqrt{2\sum_{n=1}^k b_n^{4/3}a_n^{2/3}}},
\]
\nopagebreak[3]
and using $\prod_{n=1}^k(1-\varepsilon_n^2)\geq1/2$, we have
\[
\mu_P\left(\sum_{n=1}^{k} a_n \delta_{b_n}\right)\geq  c\l(\sum\limits_{n=1}^kb_n^{4/3}a_n^{2/3}\r)^{3/4}.
\]
\end{samepage}
The same argument applies to any finite range of indices.
For general $\nu$, let $\nu_n=\nu|_{[2^{-n},2^n)}$ and round its weights
downward to obtain $\nu_{\ell}$. By the above bound and
Lemma \ref{lem: rate comparison},
\[
\mu_P(\nu)\geq\mu_P(\nu_{\ell})\geq cJ_2(\nu_n)^{3/4}.
\]
Letting $n\to\infty$ and using monotone convergence for $J_2$
gives the desired lower bound.
\end{proof}
% \begin{revision}
\begin{proof}[Proof of Lemma \ref{lem: delta_1 lb}]
Weight homogeneity and spatial scaling give
$\E L_P(a\delta_b,1)=b\lambda(\sqrt a)$.
In a square of side $s$, take $\lfloor s\rfloor$ successive unit squares
along the diagonal. Each is occupied with probability $1-e^{-1}$, and
one point from each occupied square forms an up-right chain. Hence
\[
\lambda(s)\geq(1-e^{-1})\lfloor s\rfloor
\geq(1-e^{-1})s-1.
\]
This proves the claim with $c_1=1-e^{-1}$.
\end{proof}
% \end{revision}

% \begin{revision}
\begin{proof}[Proof of Lemma \ref{smallerrectangle}]
The areas satisfy $A_{n-1,j}=l_jw_jA_{n,i}$. Put
$k=1-\varepsilon_{n-1}$. The slope constraint gives
$w_j/l_j\in[k,k^{-1}]$, and hence
\[
\frac{l_jw_j}{(l_j+w_j)^2}
\geq\frac{k}{(1+k)^2}
\geq\frac{k^2}{(1+k^2)^2}.
\]
The first inequality follows by minimizing $x/(1+x)^2$ on
$[k,k^{-1}]$; the second uses its monotonicity on $[0,1]$.
This is the asserted bound.
\end{proof}
% \end{revision}

\Needspace{7\baselineskip}
% \begin{revision}
\begin{proof}[Proof of Lemma \ref{smallerrectangles}]
The relative side lengths satisfy $\sum_jl_j=\sum_jw_j=1$.
Writing $\varepsilon=\varepsilon_{n-1}$, Lemma~\ref{smallerrectangle} gives
\[
\sum_{R_{n-1,j}\subset R_{n,i}}\sqrt{A_{n-1,j}}
\geq\frac{1-\varepsilon}{1+(1-\varepsilon)^2}
\sqrt{A_{n,i}}\sum_j(l_j+w_j)
=\frac{2(1-\varepsilon)}{1+(1-\varepsilon)^2}\sqrt{A_{n,i}}.
\]
Finally,
\[
2(1-\varepsilon)-(1-\varepsilon^2)\bigl(1+(1-\varepsilon)^2\bigr)
=\varepsilon^2(1-\varepsilon)^2\geq0,
\]
so the last coefficient is at least $1-\varepsilon^2$.
\end{proof}
% \end{revision}

\begin{proof}[Proof of Lemma \ref{parallelogram}]
    The linear transformation that takes this parallelogram to the square with opposite sides $(0,0)$ and $(\sqrt {A\alpha(\varepsilon)}, \sqrt {A\alpha(\varepsilon)})$ is area-preserving, so it preserves the Poisson point process with intensity measure $\nu$, and it sends the set of lines with slope between $m(1-\varepsilon)$ and $m/(1-\varepsilon)$ to the set of lines with positive slope.
\end{proof}

\section{Proofs of lattice model}

We prove Theorem \ref{thm:bounds for lattice model} by first proving the lower bound and then the upper bound.

\begin{proof}[Proof of Theorem \ref{thm:bounds for lattice model}]
\hspace{0.1mm}
\subsubsection*{Lower bound}
To prove the lower bound, we couple the Poisson model and the lattice model. We first suppose that $X$ is a bounded random variable. Consider the Poisson model with points $\mathcal{P}$ and the passage time {$L_P(\nu,n)$}, in the square $(0,n)\times(0,n)$ with $\nu$ corresponding to the random variable $X$. For $1\leq i,j\leq n$, define 
\[
S_{ij}=\sum\limits_{v\in\mathcal{P}\cap (i-1,i)\times (j-1,j)}w(v),\qquad M_{ij}=\max\limits_{v\in\mathcal{P}\cap (i-1,i)\times (j-1,j)}w(v),\qquad R_{ij}=S_{ij}-M_{ij}
\]

Now consider the three lattice LPP models with i.i.d vertex weights given by $S_{i,j}, M_{ij}$ and $R_{ij}$ and their LPP values from $(1,1)$ to $(n,n)$ are denoted by ${G_n}(S),{G_n}(M)$ and ${G_n}(R)$ respectively. The passage times $L_{P}(\nu,n)$ can be naturally coupled with $G_n(S),{G_n}(M),{G_n}(R)$ such that
\[L_{P}(\nu,n)\leq G_n(S)\leq G_n(M)+G_n(R).\]
Note that
$M$ is stochastically dominated by $X$.
We have that $\E G_n(M)\leq \E G_n(X)$ and {$n\E X\leq\E G_n(X)$, by evaluating any fixed up-right path.} In particular, 
\begin{equation}\label{eq: Poi-Lattice coupling bound}
    {\E L_{P}(\nu,n)\leq \E G_n(X)+\E G_n(R).}
\end{equation}
{We will show that the lattice-animal bound of
\citet[Proposition 3.4(i)]{JM04} implies $\E G_n(R)\leq\bar c n\E X$.
It follows that}
\begin{align*}
    \mu_{P}(\nu)\leq \lim_{n\rightarrow\infty}(\E G_n(M)+\E G_n(R))/n\leq \mu_{\Z^2}(\nu)+\bar c\E X\leq \bar c\mu_{\Z^2}(\nu).
\end{align*}

Then Theorem \ref{thm: bounds for poisson model} gives the desired lower bound for Theorem \ref{thm:bounds for lattice model}.

We now prove that $\E G_n(R)\leq \bar cn\E X$ for some constant $\bar c>0$. {Each up-right path has $2n-1$ vertices and is a lattice animal.
Thus \citet[Proposition 3.4(i)]{JM04}, which records the bounds of
\citet{JM02}, gives} \begin{equation}\label{eq: weak ub on LPP}
    \E G_n(R)\leq \bar cn\int_{0}^{\infty}\sqrt{\P(R>t)}dt.
\end{equation}

{Let $N$ be the number of Poisson points in a unit square.} Conditional on $N=m$, the unordered $m$ weights $X_1,X_2,...,X_m$ are independent with marginals as $X$. For $m\geq 2$, we have
\[\P(R>t|N=m)\leq \binom{m}{2}\P\l(X>\frac{t}{m-1}\r)^2.\]

As a result,
\[
\P(R>t)=\sum\limits_{m\geq 2}\P(N=m)\P(R>t|N=m).
\]
% \begin{revision}
Using $\sqrt{\sum_m u_m}\leq\sum_m\sqrt{u_m}$ and
$\P(N=m)=e^{-1}/m!$, we obtain
\[
\int_0^\infty\sqrt{\P(R>t)}\,dt
\leq \E X\sum_{m=2}^\infty(m-1)
\left(\frac{e^{-1}}{m!}\binom m2\right)^{1/2}
\leq\bar c\,\E X.
\]
The series converges by the factorial in the denominator. Consequently,
% \end{revision}
\[
\E G_n(R)\leq \bar cn\int_{0}^{\infty}\sqrt{\P(R>t)}dt\leq \bar cn\E X.
\]
{This gives the lower bound for bounded $X$. For general $X$, apply
the bound to $X\wedge k$ and let $k\to\infty$. The passage values
increase, and the rate integrals converge by monotone convergence.}
This completes the proof of the lower bound for Theorem \ref{thm:bounds for lattice model}.

\subsubsection*{Upper bound}
To prove the upper bound, we couple the lattice model and the entropy model. Let $X$ be a non-negative random variable with measure $\nu$. Let $a\geq 0$ be such that \[
a=\inf_t\{\P(X\leq t)\geq1/2\}.
\] Define \[
Y=X\mathbf{1}_{\{X>a\}}
\]. Note that
$L_{\Z^2}(X,n)\leq L_{\Z^2}(Y,n)+a(2n-1)$. We also have that \[\bar c\l(\int\limits_{0}^{\infty}\P(X>t)^{2/3}t^{1/3}dt\r)^{3/4}\geq a.\] So if we show the desired bound for $Y$, we are done.

We bound $\E L_{\Z^2}(Y,n)$ by coupling lattice model to the entropy model with measure $\eta:=2\nu|_{(a,\infty)}$ and parameter $d$. In a unit square, consider the Poisson process with measure $\eta$. {If $M$ is the maximum
weight in a unit square one can check that $M$ stochastically dominates $Y$.} We use this stochastic domination to couple $Y$ and the Poisson process independently in each square as follows.

% \begin{revision}
Consider the unit squares
$C_{ij}=(3i-2,3i-1)\times(3j-2,3j-1)$,
$1\leq i,j\leq n$, in a $3n\times3n$ grid. In each occupied square,
select a point attaining $M_{ij}$; use its centre if it is empty.
For an optimal lattice path, connect these selected locations in path
order, starting at $(0,0)$ and ending at $(3n,3n)$, and rotate as above.
An increment between consecutive locations has the form
$3e_j+\eta$, with $|\eta_\ell|\leq1$, so its time component is at least
$1/\sqrt2$ and its Dirichlet energy is bounded by a universal constant.
The initial and final increments have the same uniform energy bound.
Thus the constructed path has duration $3\sqrt2 n$, collects at least
the lattice weight, and has energy at most $\kappa n$.
Locations in empty squares may be retained as unweighted vertices;
omitting them can only decrease the energy.
% \end{revision}
We have that 
\[
L_{\Z^2}(Y,n)\leq L_D(\eta,d,(0,3\sqrt{2}n))+\kappa dn.
\]
Taking $n\rightarrow\infty$, using Theorem \ref{thm: LBCPP} and Lemma \ref{lem:entropy scaling}, and optimizing over $d$, we have that
\[
 \mu_{\Z^2}(\nu)\leq \bar c\l(\int_{0}^{\infty}\nu(x,\infty)^{2/3}x^{1/3}dx\r)^{3/4}.
\]
This completes the proof of the upper bound.
\end{proof}

\section{Proofs of higher dimension models}

\begin{proof}[Proof of Theorem \ref{thm: gen dim bounds}]
We prove the lower bounds for the Poisson model and then transfer the results to other models. For the upper bounds, we first prove it for entropy model and then transfer the results to other models.

\subsection*{Lower bounds}

We first prove the lower bounds for the Poisson model. We prove the following theorem first.
\begin{theorem}\label{thm: Poisson gen dim discrete bound}
Let \(q\ge2\), \(b_n=2^n\), and
$\nu_\ell=\sum_{n=1}^k a_n\delta_{b_n}.$ There exists \(c_q>0\) such that, for \(0<\varepsilon_n\le1/(2q)\), we have
\begin{equation}\label{eq: disc gen dim poisson bound}
    \mu_{P,q}(\nu_\ell)
\ge c_q\prod_{n=1}^k(1-(q-1)^2\varepsilon_n^2)
\sum_{n=1}^k b_na_n^{1/q}\varepsilon_n^{(q-1)/q}.
\end{equation}
\end{theorem}
The proof is similar to that of 2D case.  We pick $\{\varepsilon_n\}_{n=1}^{k}$, $\varepsilon_n\in (0,1)$ (to be decided later). We construct a path of line segments recursively. We are given a path after stage $n$ (chronologically step indexed by $k$ is performed first and then $k-1$ and so on), which goes through vertices of weight at least $b_{n}$. This set of vertices is denoted $E_n$, and is ordered naturally so that one obtains an up-right path. Initially, i.e. after stage $k+1$, the set of vertices is $E_{k+1} = \{ \mathbf{0}, (t,\ldots,t)\}$.

At stage $n$, we run Poisson last passage on $\nu = a_{n} \delta_{b_{n}}$, between every consecutive pair of points in $E_{n+1}$, with an additional constraint. {Let the parent box have side lengths \(L_1,\ldots,L_q\).
For an increment \(\Delta x\) inside the box, define its relative coordinates by}
\[
r_\ell=\frac{\Delta x_\ell}{L_\ell},
\qquad
\bar r=\frac1q\sum_{\ell=1}^qr_\ell.
\]
The constraint is,
\begin{equation}\label{eq: gen dim constraint}
{r_\ell\ge(1-\varepsilon)\bar r
\quad\text{for every }\ell.}    
\end{equation}
 If the points picked up by this last passage are denoted by $P_{n}$, then $E_{n} = P_{n} \cup E_{n+1}$, and then we proceed to step $n-1$. We first compute a lower bound on $\E L(\sum_{i=1}^{k} a_i \delta_{b_i},t)$. 
% \begin{revision}
The points in $E_i$ are totally ordered by the coordinates. Write
$p_s=(x^{(1)}_{i,s},\ldots,x^{(q)}_{i,s})$ for its $(s+1)$st point.
The first and last points remain in $E_{i-1}$. Define the box
\[
R_{i,s}=\prod_{\ell=1}^q
[x^{(\ell)}_{i,s-1},x^{(\ell)}_{i,s}],
\]
with volume $A_{i,s}$, for $1\leq s\leq|E_i|-1$.
% \end{revision}

Since $E_{k+1}= \{ \mathbf{0}, (t,\ldots,t)\}$, at stage $k$, between these two points in $E_{k+1}$, we perform Poisson LPP on $\nu=a_k\delta_{b_k}$ with above mentioned constraint. We use the following lemmas, the proofs of which are postponed.

    \begin{lemma}\label{lem: delta_1 lb gen dim}
There exists $c_1>0$ so that 
$$
\E L_{P,q}(a \delta_b,1) \geq b(c_1 a^{1/q} - 1)
$$
for all $a,b>0$.
\end{lemma}

\begin{lemma} \label{smallerrectangle gen dim}
    {Consider a child box $R_{n-1,j}\subset R_{n,i}$.
The quantities $x_{n-1,j}^{(\ell)}$ are vertex coordinates; define the relative side lengths by}
\[
l_{j,\ell}=\frac{x_{n-1,j}^{(\ell)}-x_{n-1,j-1}^{(\ell)}}
 {x_{n,i}^{(\ell)}-x_{n,i-1}^{(\ell)}},
\qquad s_j=\frac1q\sum_{\ell=1}^q l_{j,\ell}.
\]
Suppose,
\[
l_{j,\ell}\ge
(1-\varepsilon_{n-1})s_j\qquad\text{for every }\ell,
\]
then,
    \begin{equation*}
{
A_{n-1,j}\ge
s_j^q
{(1-\varepsilon_{n-1})^{q-1}
\bigl(1+(q-1)\varepsilon_{n-1}\bigr)}
A_{n,i}.
}
    \end{equation*}
\end{lemma}

\begin{lemma} \label{smallerrectangles gen dim}
    {For any parent box $R_{n,i}$, we have}
    \begin{equation*}
        \sum_{R_{n-1,j} \subset R_{n,i}} {A_{n-1,j}^{1/q}} \geq \left(1-(q-1)^2{\varepsilon_{n-1}^2}\right){A_{n,i}^{1/q}}.
    \end{equation*}
\end{lemma}

\begin{lemma}\label{parallelogram gen dim}
{Define $l_{\varepsilon}(A\nu)$ to be the last passage time in a box
of volume $A$ with intensity measure $\nu$, subject to
\eqref{eq: gen dim constraint} in relative coordinates.} Then
\begin{equation*}
    l_{\varepsilon}(A\nu) \stackrel{d}{=} L_{P,q}(\nu,(\varepsilon^{q-1}A)^{1/q}).
\end{equation*}
\end{lemma}
\begin{proof}[Proof of Theorem \ref{thm: Poisson gen dim discrete bound}]

Let $\mathcal F_n$ be the $\sigma$-field generated by all points of weight at least $b_n$. For $i\in[|E_n|-1]$, let $X_{n-1,i}$ be the number of points $P_{n-1}\cap R_{n,i}$. Let \[
\lambda_q(s)=\mathbb E L_{P,q}(\delta_1,s).
\]
We have, by Lemma \ref{parallelogram gen dim},
$$
\E[X_{n-1,i}|\mathcal F_n]=\lambda_q\l(\l(a_{n-1}A_{n,i}\varepsilon_{n-1}^{q-1}\r)^{1/q}\r),
$$
giving the bound
\begin{align}\label{eq: lb for number of points gen dim}
   \E[X_{n-1,i}|\mathcal F_n]\ge c_1a_{n-1}^{1/q}A_{n,i}^{1/q}\varepsilon_{n-1}^{(q-1)/q}-1, 
\end{align}
by Lemma \ref{lem: delta_1 lb gen dim}. If $L_{n-1}$ is the contribution to  $L_q(\sum_{n=1}^{k} a_n \delta_{b_n},t))$ just from the points with weight $b_{n-1}$, then using \eqref{eq: lb for number of points gen dim} and Lemma \ref{smallerrectangles gen dim},
\begin{align*}
    \E[L_{n-1}|\mathcal{F}_n]&=b_{n-1}\sum_{i=1}^{|E_n|-1}\E[X_{n-1,i}|\mathcal{F}_{n}]\\
    &\geq
c_0t\,b_{n-1}a_{n-1}^{1/q}
\varepsilon_{n-1}^{(q-1)/q}
\prod_{j=n}^k(1-(q-1)^2\varepsilon_j^2)-b_{n-1}(|E_n|-1).
\end{align*}
% \begin{revision}
Put $W=\sum_{n=1}^kL_n$. The same geometric-weight bookkeeping
as in the planar proof gives
$\sum_{n=1}^k b_n(|E_{n+1}|-1)\leq2b_k+W$.
Taking expectations and summing therefore yields
\[
2\E W\geq c_0t\prod_{n=1}^k(1-(q-1)^2\varepsilon_n^2)
\sum_{n=1}^k b_na_n^{1/q}\varepsilon_n^{(q-1)/q}-2b_k.
\]
Since $L_{P,q}(\nu_\ell,t)\geq W$, we have the finite-size estimate
\begin{equation}\label{eq: gen dim poisson finite size estimate}
\E L_{P,q}(\nu_\ell,t)
\geq \frac{c_0t}{2}\prod_{n=1}^k(1-(q-1)^2\varepsilon_n^2)
\sum_{n=1}^k b_na_n^{1/q}\varepsilon_n^{(q-1)/q}-b_k.
\end{equation}
Division by $t$ and passage to the limit complete the proof of
Theorem~\ref{thm: Poisson gen dim discrete bound}.
% \end{revision}
\end{proof}
The same optimization over $\varepsilon_n$ as in the 2D case, followed by truncation and Lemma \ref{lem: rate comparison}, gives us
\begin{equation}\label{eq: Poi gen dim lb}
    \mu_{P,q}(\nu)\geq c_{q,P}\l(\int_{0}^{\infty}\nu(x,\infty)^{2/(q+1)}x^{(q-1)/(q+1)}dx\r)^{{(q+1)/(2q)}}
\end{equation}

% \begin{revision}
We now transfer this lower bound to the entropy and lattice models.
Rotate the time direction to $q^{-1/2}\mathbf1$, so that the points
$t e_i$ have time coordinate $t/\sqrt q$, while $t\mathbf1$ has time
coordinate $t\sqrt q$ and zero spatial coordinates. If $(\Delta z,\Delta u)$
is a rotated up-right increment, then
\[
\|\Delta z\|^2\leq(q-1)(\Delta u)^2.
\]
Consequently, the rotated path has Dirichlet energy at most
$(q-1)\sqrt q\,t$, and penalty at most $d(q-1)\sqrt q\,t$. Thus
% \end{revision}
\begin{equation}\label{eq: Poisson-entropy lb gen dim}
    L_{D,q}(\nu,d,\sqrt{q}t)\geq L_{P,q}(\nu,t)-d(q-1)\sqrt{q}t. 
\end{equation}
Dividing by $t$, taking $t\rightarrow\infty$, and using \eqref{eq: Poi gen dim lb} and Lemma \ref{lem:entropy scaling}, optimization over $d$ gives 
\begin{equation}\label{eq: entropy lb gen dim}
\mu_{D,q}(\nu,1)
\ge
c_{q,D}\int_{0}^{\infty}\nu(x,\infty)^{2/(q+1)}x^{(q-1)/(q+1)}dx
\end{equation}

% \begin{revision}
For the lattice lower bound, let $S_v$ be the sum of Poisson marks
in a unit cube, $M_v$ its largest mark, and $M_{v,q}$ the sum of its
largest $q-1$ marks, using zeros if fewer marks occur. Put
$R_{v,q}=S_v-M_{v,q}$. Then $M_{v,q}\leq(q-1)M_v$ and
$M_v$ is stochastically dominated by $X$.
Conditional on the number $N=m\geq q$ of marks in the cube,
\[
\P(R_{v,q}>t\mid N=m)
\leq\binom mq\P\left(X>\frac{t}{m-q+1}\right)^q.
\]
Indeed, if the sum remaining after removal of the largest $q-1$ marks
exceeds $t$, the $q$th largest mark exceeds $t/(m-q+1)$.
Since $N$ is Poisson with mean one,
\[
\int_0^\infty\P(R_{v,q}>t)^{1/q}\,dt
\leq\E X\sum_{m=q}^\infty(m-q+1)
\left(\frac{e^{-1}}{m!}\binom mq\right)^{1/q}
\leq c_q\E X.
\]
The series converges. The lattice-animal bound in
\citet[Proposition 3.4(i)]{JM04} now gives the finite-size comparison
\begin{equation}\label{eq: poisson lattice finite q}
\E L_{P,q}(\nu,n)
\leq(q-1)\E L_{\Z^q}(\nu,n)+c_q n\E X.
\end{equation}
Every fixed up-right path gives $\E L_{\Z^q}(\nu,n)\geq n\E X$.
Divide by $n$, take limits, and use \eqref{eq: Poi gen dim lb} to obtain
\[
\mu_{\Z^q}(\nu)\geq c_q J_q(\nu)^{(q+1)/(2q)}.
\]
As in the planar case, first argue for bounded $X$ and then use
$X\wedge k$ and monotone convergence.
% \end{revision}
This completes the proof of lower bound for all the three models. We now prove the upper bounds.

\subsection*{Upper bounds}
We first prove the upper bounds for the entropy model with parameter $d=1$.
Defining $W_{a,w,\theta}(\mathbf{x},\mathbf{y};t)$ similar to the 2D case and the same proof of Lemma \ref{lem:insert}, gives us that for every $\mathbf{x},\mathbf{y}\in\R^{q-1}$ and $a,w,\theta,t>0$,
\begin{equation*}
\E W_{a,w,\theta}(\mathbf{x},\mathbf{y};t)
\le \exp\left\{
\bar c_q a^{2/(q+1)}
\theta ^{-(q-1)/(q+1)}
e^{2\theta w/(q+1)}t
\right\},
\end{equation*}
for some universal constant $\bar c_q$. Defining $W_{\theta}(\sigma)$ as in 2D case gives us analog of Corollary \ref{cor:polygon} as,
\[
\mathbb E W_\theta(\sigma)
\le
\exp\left\{
\bar c_q a^{2/(q+1)}
\theta ^{-(q-1)/(q+1)}
e^{2\theta w/(q+1)}t
\right\}.
\]
Using the above bounds, repeating the same arguments of the 2D case, and applying Lemma \ref{lem: rate comparison}, we get that 
\begin{equation}\label{eq: entropy ub gen dim}
\mu_{D,q}(\nu,1)
\le
\bar c_{q,D}\int_{0}^{\infty}\nu(x,\infty)^{2/(q+1)}x^{(q-1)/(q+1)}dx.
\end{equation}
The same arguments of 2D case of coupling entropy model with Poisson model in higher dimensions, in particular using \eqref{eq: Poisson-entropy lb gen dim}, Lemma \ref{lem:entropy scaling} and optimization over $d$, gives us the desired upper bound for the Poisson model.

% \begin{revision}
For the lattice upper bound, repeat the median truncation and
Poisson-maximum coupling from the planar case. Put $k=2q-1>q$ and use
the separated unit cubes
\[
C_v=\prod_{\ell=1}^q(kv_\ell-2,kv_\ell-1),
\qquad v\in\{1,\ldots,n\}^q,
\]
in a grid of side length $kn$, containing $(kn)^q$ unit cubes.
Select a maximum-weight point in each occupied cube and its centre
in an empty cube. Consecutive locations along a lattice path differ by
$k e_j+\eta$, where $|\eta_\ell|\leq1$. After rotation,
\[
\Delta u\geq\frac{k-q}{\sqrt q}>0,
\qquad \|\Delta z\|\leq k+\sqrt q,
\]
so each segment has Dirichlet energy bounded by a constant depending
only on $q$. The initial and final connections to the lower and upper
corners of $[0,kn]^q$ obey the same bound. There are $O_q(n)$ segments,
so the total penalty is at most $c_qdn$, and the duration is
$k\sqrt q\,n$. Applying \eqref{eq: entropy ub gen dim} and
Lemma~\ref{lem:entropy scaling}, then optimizing over $d$, gives the
lattice upper bound. This completes the proof of
Theorem~\ref{thm: gen dim bounds}.
% \end{revision}
\end{proof}

\begin{proof}[Proof of Lemma \ref{lem: delta_1 lb gen dim}]
    Taking $\lfloor a^{1/q}\rfloor$ many successive cubes of side length $a^{-1/q}$, each of which is non-empty with probability $1-e^{-1}$, we get 
    \[\E L_{P,q}(a \delta_b,1) \geq b(1-e^{-1}) {\lfloor a^{1/q}\rfloor}.\]
    This implies the desired bound. 
\end{proof}

\begin{proof}[Proof of Lemma \ref{smallerrectangle gen dim}]
     The volumes are related by
     \[
A_{n-1,j}=A_{n,i}\prod_{\ell=1}^q l_{j,\ell}.
\]
{Keeping $qs_j$ fixed, the product is minimized when $q-1$
coordinates equal $(1-\varepsilon_{n-1})s_j$ and the remaining coordinate
equals $(1+(q-1)\varepsilon_{n-1})s_j$. Indeed, with all other coordinates
fixed, the product is a concave quadratic in either of two coordinates
whose sum is fixed, so one of them may be moved to its lower bound.
Iterating and substituting these values gives the stated inequality.}
   
\end{proof}

% \begin{revision}
\begin{proof}[Proof of Lemma \ref{smallerrectangles gen dim}]
Write $\varepsilon=\varepsilon_{n-1}$ and
$x=(q-1)\varepsilon$. In the range $0<\varepsilon\leq1/(2q)$,
Bernoulli's inequality gives
\[
(1-\varepsilon)^{q-1}(1+(q-1)\varepsilon)
\geq(1-x)(1+x)=1-x^2.
\]
As $0\leq1-x^2\leq1$, its $q$th root is at least $1-x^2$.
Moreover, the child increments sum to the parent increment in every
coordinate, so
\[
\sum_j s_j=\frac1q\sum_{\ell=1}^q\sum_j l_{j,\ell}=1.
\]
Take $q$th roots in Lemma~\ref{smallerrectangle gen dim} and sum over $j$.
This yields the claimed factor $1-(q-1)^2\varepsilon^2$.
\end{proof}
% \end{revision}

% \begin{revision}
\begin{proof}[Proof of Lemma \ref{parallelogram gen dim}]
Normalize the parent box to $[0,1]^q$ by dividing each coordinate
by the corresponding side length, and consider
\[
B=\varepsilon I+\frac{1-\varepsilon}{q}\mathbf1\mathbf1^T.
\]
The map preserves coordinate averages, fixes $\mathbf1$, and has
determinant $\varepsilon^{q-1}$. For $z=Bv$,
\[
v=\frac{z-(1-\varepsilon)\bar z\,\mathbf1}{\varepsilon}.
\]
Thus $v\geq0$ if and only if $z$ satisfies
\eqref{eq: gen dim constraint}. A constrained path from $0$ to
$\mathbf1$ is therefore exactly the image under $B$ of an up-right
path in $[0,1]^q$. The marked intensity in the latter coordinates is
$A\varepsilon^{q-1}\nu$. An isotropic rescaling turns this into intensity
$\nu$ in a cube of side $(A\varepsilon^{q-1})^{1/q}$, proving the claim.
\end{proof}
% \end{revision}

\Needspace{7\baselineskip}
\section{Proofs of applications}

\begin{proof}[Proof of Theorem \ref{thm: critical heavy tailed LPP bound}]

We first prove the upper bound and then the lower bound. We prove for $q=2$. The same argument works for any $q\geq 2$. Take $Y=X\wedge  n$. Let $G_n(X), G_n(Y)$ be the passage times with i.i.d vertex weights $X$ and $Y$, from $(1,1)\rightarrow(n,n).$ We have,
\[
G_n(X)\leq G_n(Y)+\sum \limits_{v\in\{1,\dots,n\}^2}
(X_v-n)_+.\]

By upper bound in Theorem \ref{thm:bounds for lattice model} and {superadditivity of $\E G_n(Y)$ with the law of $Y$ held fixed} there exists a constant $\bar c$, independent of $n$, such that \[\E G_n(Y)\leq \bar cn\l(\int_{0}^n\P(X>t)^{2/3}t^{1/3}dt\r)^{3/4}\leq \bar cn\l(\log n\r)^{3/4}.\] Also,
\[
\E\sum_{v\in\{1,\dots,n\}^2}(X_v-n)_+
\leq n^2\int_n^\infty\P(X>t){\,dt}\leq\bar c n.
\]
This gives the upper bound in Theorem \ref{thm: critical heavy tailed LPP bound}. We now prove the lower bound.

Let {$a_j=\P\l(2^j<X<2^{j+1}\r)$} and let $b_j=2^j$ and $k=\lfloor \log_2n\rfloor.$ Choose $j_0$ such that $0<c_1\leq b_j\sqrt{a_j}\leq c_2<\infty$ for all $j\geq j_0.$ There are $m=k-j_0+1$ levels and let $\eta=\sum_j a_j\delta_{b_j}$. 
% \begin{revision}
Take $\varepsilon_j=(4\sqrt m)^{-1}$ for $j_0\leq j\leq k$. Then
\[
\prod_{j=j_0}^k(1-\varepsilon_j^2)\geq\frac{15}{16},
\qquad
\sum_{j=j_0}^k b_j\sqrt{a_j\varepsilon_j}\asymp m^{3/4}.
\]
% \end{revision}
Applying $q=2$ estimate in \eqref{eq: gen dim poisson finite size estimate}, we get
\[
\E L_P(\eta,n)\geq cn\l(\log n\r)^{3/4}-\bar cn.
\]
 Since $\E X<\infty$, applying \eqref{eq: Poi-Lattice coupling bound} (it can be shown that similar estimate holds for higher dimensions), we have the lower bound of Theorem \ref{thm: critical heavy tailed LPP bound}.    
\end{proof}

\paragraph{AI assistance.}
The authors developed the idea of using entropy-controlled model, Poisson model instead of the lattice model. The authors also derived the exact tail integral condition, the proof of the lower bound and the upper bound proof strategy. A simpler upper bound proof was found together with ChatGPT Astra 6.1, which also helped with final editing and error checking. Claude Opus 5.5 was used to create figures.  

\paragraph{Acknowledgments.}
The authors were supported by B.V.'s NSERC Discovery Accelerator grant.

% Keep the end matter together without changing line spacing.
\setlength{\bibsep}{4pt}
\bibliography{shape}
\bibliographystyle{plainnat-linked}

\par\vspace{2\baselineskip}
\noindent\begin{minipage}{\textwidth}
\normalsize
\noindent
Jnaneshwar Baslingker. Department of Mathematics, University of Toronto,
40 St. George St., Toronto, Ontario, M5S 2E4, Canada.
\href{mailto:j.baslingker@utoronto.ca}{\texttt{j.baslingker@utoronto.ca}}
\par\medskip
\noindent
Daniel Spivak. Department of Mathematics, University of Toronto,
40 St. George St., Toronto, Ontario, M5S 2E4, Canada.
\href{mailto:daniel.spivak@mail.utoronto.ca}{\texttt{daniel.spivak@mail.utoronto.ca}}
\par\medskip
\noindent
B\'alint Vir\'ag. Department of Mathematics, University of Toronto,
40 St. George St., Toronto, Ontario, M5S 2E4, Canada.
\href{mailto:balint@math.toronto.edu}{\texttt{balint@math.toronto.edu}}
\par
\end{minipage}

\end{document}